\documentclass[onecolumn]{elsart3p}

\journal{arXiv.org}

\newtheorem{theorem}{Theorem}[section]
\newtheorem{lemma}[theorem]{Lemma}
\newtheorem{remark}[theorem]{Remark}
\newtheorem{example}[theorem]{Example}

\newfont{\cirilrm}{wncyr10 scaled 1000}
\newfont{\cirilit}{wncyi10 scaled 1000}

\newcommand{\ch}{\symbol{'161}}

\newcommand{\ja}{\symbol{'37}}

\newcommand{\mek}{\symbol{'176}}

\newcommand{\ijj}{\symbol{'32}}

\newcommand{\ii}{\symbol{'171}}

\begin{document}

\begin{frontmatter}
\title{\bf Formulae of Partial Reduction for Linear Systems \\ of First Order Operator Equations}

\author[etf]{\bf Branko Male\v sevi\' c},
\ead{malesevic@etf.rs}
\author[matf]{\bf Dragana Todori\' c},
\ead{draganat@matf.bg.ac.rs}
\author[etf]{\bf Ivana Jovovi\' c},
\ead{ivana@etf.rs}
\author[matf]{\bf Sonja Telebakovi\' c}
\ead{sonjat@matf.bg.ac.rs}

\vspace{10pt}
\address[etf]{Faculty of Electrical Engineering, University of  Belgrade,
Bulevar kralja Aleksandra 73, 11000 Belgrade, Serbia}
\address[matf]{Faculty of Mathematics, University of  Belgrade,
Studentski trg 16, 11000 Belgrade, Serbia}

\begin{abstract}
This paper deals with reduction of non-homogeneous linear systems of first order operator
equations with constant coefficients. An equivalent reduced system, consisting of higher order linear operator
equations having only one variable and first order linear operator equations in two variables, is obtained by using
the rational canonical form.
\end{abstract}

\begin{keyword}
Linear system of first order operator equations with constant coefficients,
$n^{th}$ order linear operator equation with constant coefficients,
sum of principal minors, the rational canonical form, the characteristic polynomial
\end{keyword}
\end{frontmatter}

\section{Introduction}

Let $K$ be a field, $V$ a vector space over $K$ and let $A:V\longrightarrow V$ be a
linear operator.
{\boldmath \bf Linear system of first order $A$-operator equations with constant coefficients
in unknowns $x_{i}$, $1 \!\leq\! i \!\leq\! n$,}~is$:$
\begin{equation}
\label{B-Sistem-1}
\left\{
\begin{array}{rcl}
A(x_1) & \!\!=\!\!      & b_{11} x_1 + b_{12} x_2 + \ldots + b_{1n} x_n + \varphi_{1} \\
A(x_2) & \!\!=\!\!      & b_{21} x_1 + b_{22} x_2 + \ldots + b_{2n} x_n + \varphi_{2} \\
       & \!\!\vdots\!\! &                                                             \\
A(x_n) & \!\!=\!\!      & b_{n1} x_1 + b_{n2} x_2 + \ldots +
b_{nn} x_n + \varphi_{n}
\end{array}
\right\}
\end{equation}
for $b_{ij}\!\in\!K$ and $\varphi_{i}\!\in\!V$. We say that $B\!=\![b_{ij}]_{i, j=1}^n\!\in\!K^{n \times n}$
is {\bf the system matrix}, and $\vec{\varphi}\!=\![\varphi_{1}\;\ldots\;\varphi_{n}]^T\!\in\!V^{n\times 1}$
is {\bf the free column}. If $\vec{\varphi} \!=\! \vec{0}$
system (\ref{B-Sistem-1}) is called homogeneous.
Otherwise it is called non-homogeneous.

Let $\vec{x}\!=\![x_{1}\;\ldots\;x_{n}]^T$ be a column of unknowns and
$\vec{A} : V^{n\times 1}\longrightarrow V^{n\times 1}$ be a vector operator
defined componentwise $\vec{A}(\vec{x})\!=\![A(x_1)\;\ldots\;A(x_n)]^T.$
Then system (\ref{B-Sistem-1}) can be written in the following vector form$:$
\begin{equation}
\label{B-Sistem-2}
\hspace{10pt}
\vec{A}(\vec{x})=B\vec{x}+\vec{\varphi}.
\end{equation}
A solution of vector equation (\ref{B-Sistem-2}) is any column $\vec{v}\!\in\! V^{n \times 1}$ which satisfies (\ref{B-Sistem-2}).

\break

Powers of operator $A$ are defined as usual $A^{i}\!=\!A^{i-1}\!\circ A$
assuming that $A^{0}\!=\!I:V \longrightarrow V$ is the identity operator.
By {\bf \boldmath $n^{th}$ order linear $A$-operator equation with constant coefficients, in unknown $x$},
we mean$:$
\begin{equation}
\label{Eq_5}
A^{n}(x) + b_{1} A^{n-1}(x) + \ldots + b_{n} I(x) =\varphi,
\end{equation}
where $b_1, \ldots, b_n \!\in\! K$ are coefficients and
$\varphi \!\in\! V$.
If $\varphi\!=\!0$ equation (\ref{Eq_5})
is called homogeneous.
Otherwise it is non-homogeneous.
A solution of equation ({\ref{Eq_5}) is any vector $v\!\in\!V$ which satisfies (\ref{Eq_5}).

\section{Notation}

Let $M \!\in\! K^{n \times n}$  be an $n$-square matrix. We denote by$:$

\vspace{5pt}
\begin{description}
\item[($\mbox{\large $\delta$}_{k}$)]
$\mbox{\large $\delta$}_{k} \!=\!\mbox{\large $\delta$}_{k}(M)$
the sum of its principal minors of order $k$ ($1 \leq k \leq n$),
\item[($\mbox{\large $\delta$}^{i}_{k}$)]
$\mbox{\large $\delta$}^{i}_{k}\!=\!\mbox{\large $\delta$}^{i}_{k}(M)$
the sum of its principal minors of order $k$ containing $i^{th}$ column
$(1 \leq i, k \leq n)$.
\end{description}

\vspace{5pt}
Let $v_{1}, \ldots, v_{n}$ be elements of $K$. We write
$M^{i} (v_{1}, \ldots, v_{n})\!\in\! K^{n \times n}$ for
the matrix obtained by substituting $\vec{v}\!=\![ v_{1}\;\ldots\; v_{n} ]^{T}$
in place of $i^{th}$ column of $M$. Furthermore, it is convenient to use$:$
$$
\hspace{-193pt}
\mbox{\large $\delta$}^{i}_{k}(M; \vec{v})
=\mbox{\large $\delta$}^{i}_{k}(M; v_{1}, \ldots, v_{n})
=\mbox{\large $\delta$}^{i}_{k}(M^i (v_{1}, \ldots, v_{n})).
$$

{\bf The characteristic polynomial} $\Delta_{B}(\lambda)$
of the matrix $B \!\in\!K^{n \times n}$ has the form
$$
\hspace{-168pt}
\Delta_{B}(\lambda)=\det(\lambda I - B)=
\lambda^{n} + d_{1} \lambda^{n-1} + \ldots + d_{n-1} \lambda +d_{n},
$$
where $d_{k} = (-1)^{k}\mbox{\large$\delta$}_{k}(B)$, $1 \leq k \leq n$; see {[{\bf 2}, p.$\,$78]}.

\vspace{5pt}
{\bf The companion matrix of a monic polynomial}
$p(\lambda) \!=\!\lambda^{n} + d_{1} \lambda^{n-1} + \ldots + d_{n-1} \lambda +d_{n}$
is the matrix
$$
\hspace{-19pt}
C_p \!=\!
\left[
\begin{array}{ccccc}
0      & 1        & \ldots & 0      & 0          \\
0      & 0        & \ldots & 0      & 0          \\
\vdots & \vdots   &        & \vdots & \vdots     \\
0      & 0        & \ldots & 0      & 1          \\
-d_{n} & -d_{n-1} & \ldots & -d_{2}  & -d_{1}
\end{array}
\right].
$$

The matrix $C$ is in {\bf the rational canonical form} if it is block diagonal,
$
C = \mbox{\rm diag}{\big (}C_{p_{1}}, C_{p_{2}}, \ldots, C_{p_{k}}{\big )}
$,
where $C_{p_{i}}$ is the companion matrix of a monic polynomial
$p_{i}(\lambda)\!=\!\lambda^{n_{i}} + d_{i,1}\lambda^{n_{i}-1} +\ldots + d_{i,n_{i}-1} \lambda + d_{i,n_{i}}$
of degree $n_{i} \!=\!\deg(p_{i}) \geq 1$,
and $p_{i} \;|\; p_{i+1}$ for all $i$ such that $1 \!\leq\! i \!<\! k$.

\section{Main Results}

The reduction of system (\ref{B-Sistem-2}) to a {\bf partially reduced system} will be divided into two steps.

Firstly, by using some basic properties of the rational canonical form,
the initial system may be transformed into a proper system for further study.

\begin{theorem}{\rm [{\sl On rational canonical form} {\bf [1{\sl$,$ p.}\,$397$]}]}
\label{Teorema_RKF}
Any square matrix $B\!\in\!K^{n \times n}$ is similar to a unique matrix $C$ in the rational canonical form.
The matrix $C$ is called the rational canonical form of the matrix~$B$.
\end{theorem}

\vspace{5pt}
\begin{lemma}
\label{Lema1}
Let $C$ be the rational canonical form of the matrix $B$, i.e. there exists a regular matrix $P$ such that $C=P^{-1} B P$.
Then the system given in the form
{\rm (\ref{B-Sistem-2})}, $\vec{A}(\vec{x}) = B \vec{x} + \vec{\varphi}$,
can be reduced to
the system
\begin{equation}
\label{C-Sistem-0}
\vec{A}(\vec{y}) = C \vec{y} + \vec{\psi},
\end{equation}
where $\vec{\psi} = P^{-1} \vec{\varphi}$ is its free column
and $\vec{y} = P^{-1} \vec{x}$ is a column of the unknowns.
\end{lemma}

\pagebreak

The following lemma will be very useful in the second step of our reduction process.

\medskip

\begin{lemma}
\label{Lema3}
If the matrix $M$ has the form$:$
$$
 \!\!\!\! M = \left [
\begin{array}{ccccc}
b_{1}     & 1       & 0       & \ldots & 0      \\[0.5 ex]
b_{2}     & 0       & 1       & \ldots & 0      \\[0.5 ex]
\vdots    & \vdots  & \vdots  &        & \vdots \\[0.5 ex]
b_{n-1}   & 0       & 0       & \ldots & 1      \\[0.5 ex]
b_{n}     & a_{n-1} & a_{n-2} & \ldots & a_{1}
\end{array}
\right ]
$$
then it follows
$$
\hspace{68pt}
\mbox{\large $\delta$}^{1}_{k}(M)=
(-1)^{k} {\Big (}
\displaystyle\sum\limits_{j=1}^{k-1}{b_{j} a_{k\!-\!j}} \!-\!b_{k}
{\Big )}  \qquad (1 \leq k \leq n).
$$
\end{lemma}

\noindent
{\bf Proof.}
Observe that the function $\mbox{\large$\delta$}^{1}_{k}$
is linear with respect to the first column of its matrix argument.
The first column $[b_{1} \; \ldots \; b_{n}]^T $ of the matrix M is equal to $\sum_{j=1}^{n}{b_{j} \vec{e}_{j}}$,
where $\vec{e}_{j}$ denotes $j^{th}$ column of the identity matrix of order $n$.
Therefore
$\mbox{\large $\delta$}^{1}_{k}(M)=\sum_{j=1}^{n}{b_{j}\mbox{\large $\delta$}^{1}_{k}(M; \vec{e}_{j})}$.
It remains to prove that
$$
\mbox{\large $\delta$}^{1}_{k}(M; \vec{e}_{j})=
\left\{
\begin{array}{ll}
0,              & k<j;    \\[1.0 ex]
-(-1)^{k},        & k=j;  \\[1.0 ex]
(-1)^k a_{k-j}, & k>j.
\end{array}
\right.
$$
Let $\alpha_{k, j}$ be an arbitrary nonzero minor in the sum
$\mbox{\large $\delta$}^{1}_{k}(M; \vec{e}_{j})$.
Since it is a nonzero principal minor which contains $1$ from the first column,
it must include $1$ from the first row. Hence, minor  $\alpha_{k, j}$ must also
include the element $1$ from $2^{nd}$ row and $3^{rd}$ column.
In the same manner we conclude that our minor
contains $1$ from $i^{th}$ row and $(i\!+\!1)^{th}$ column for $i\in \{2, \ldots, j-1 \}$.
Accordingly, minor $\alpha_{k, j}$ must contain the first $j$ columns and rows.

\smallskip
\noindent
For $k\!<\!j$ there is no such nonzero principal minor.

\smallskip
\noindent
If $k\!=\!j$, then the only nonzero minor in the sum
$\mbox{\large $\delta$}^{1}_{k}(M; \vec{e}_{j})$ is the minor
formed from the first $j$ columns and rows. Expanding it
by the first column, we easily conclude that its value is
$(-1)^{k-1}$.

\smallskip
\noindent
For $k\!>\!j$ the same reasoning implies that
$\alpha_{k, j}$ must contain all columns and
rows from $1^{st}$ to $j^{th}$.
Now, let $l\!+\!1$ be the smallest index greater than $j$ such
that $\alpha_{k, j}$ contains $(l\!+\!1)^{th}$ column.
Consequently, the minor must include element $a_{n-l}$.
We apply zig-zag argument again to deduce that
our minor has a block formed from the last $n\!-\!l$ columns and rows.
So, the only minor having the desired properties is the one obtained
by deleting all columns and rows between $(j\!+\!1)^{th}$ and $l^{th}$.
It follows that
$k\!=\!j+n-l$ and $\alpha_{k,j}=(-1)^{j-1}\cdot(-1)^{n-l-1}a_{n-l}=(-1)^{k}a_{k-j}. \; \Box$

\vspace{5pt}
\begin{remark}{\rm The matrix $M$ is the transpose of the doubly companion matrix introduced in {\bf [4]}, see also~\mbox{\bf [5-7]}.}
\end{remark}

\medskip

We will now prove the main theorem dealing with a special case when the rational canonical form
of the~system matrix  has only one block.

\medskip

\begin{theorem}
\label{Teorema_parc_red_1}
Assume that the rational canonical form of the matrix $B$ is exactly the companion matrix $C$
of the characteristic polynomial
$\Delta_{B}(\lambda) \!=\! \lambda^{n} + d_{1} \lambda^{n-1} +
\ldots + d_{n-1} \lambda + d_{n} \!=\! \Delta_{C}(\lambda)$.
Then linear system of first order $A$-operator equations {\rm (\ref{B-Sistem-1})}
can be reduced to the partially reduced system$:$

\break

\noindent
\begin{equation}
\label{PRS 1}
\!\!(\mbox{$\cal R$}_{C})\!\! \quad
\left\{
\begin{array}{ccl}
\Delta_{C}(A)(y_{1})   &\!\!=\!\!& \;
\displaystyle\sum\limits_{k=1}^{n}{\!(-1)^{k+1}
\mbox{\large $\delta$}^{1}_{k}(C; A^{n-k}(\psi_{1}), \ldots, A^{n-k}(\psi_{n}))\!}     \\[0.0 ex]
y_{2}                  &\!\!=\!\!& \; A(y_{1}) - \psi_{1}                              \\[0.0 ex]
y_{3}                  &\!\!=\!\!& \; A(y_{2}) - \psi_{2}                              \\[0.0 ex]
                       & \vdots  &                                                     \\[0.0 ex]
y_{n}                  &\!\!=\!\!& \; A(y_{n-1}) - \psi_{n-1}
\end{array}
\right\},
\end{equation}
where the columns  $\vec{y}\!=\![y_{1}\;\ldots\;y_{n}]^{T}$ and
$\vec{\psi} \!=\![\psi_{1}\;\ldots\;\psi_{n}]^{T}$ are determined by
$\vec{y} \!=\! P^{-1} \vec{x}$ i $\vec{\psi} \!=\! P^{-1}\vec{\varphi}$
for a~regular matrix $P$ such that $C \!=\!P^{-1} B P$.
\end{theorem}

\smallskip

\noindent
{\bf Proof.}
Applying Lemma \ref{Lema1}, system (\ref{B-Sistem-2}) is equivalent to system (\ref{C-Sistem-0}),
which can be presented in the following form$:$

\vspace{-10pt}
\begin{equation}
\label{C-Sistem-1}
\phantom{\!\!(\mbox{$\cal R$}_{C})\!\! \quad }
\left\{
\begin{array}{ccl}
A(y_{1})   &\!\!=\!\!& \; y_{2} + \psi_{1}                                            \\[0.0 ex]
A(y_{2})   &\!\!=\!\!& \; y_{3} + \psi_{2}                                            \\[-0.5 ex]
           & \vdots  &                                                                \\[-0.5 ex]
A(y_{n-1}) &\!\!=\!\!& \; y_{n} + \psi_{n-1}                                          \\[0.0 ex]
A(y_{n})   &\!\!=\!\!& -d_{n} y_{1} - d_{n-1} y_{2} - \ldots -
d_{1} y_{n} + \psi_{n}
\end{array}
\right\}.
\end{equation}
We now consider the linear operator:

\vspace{-10pt}
$$
\hspace{-222pt}
\Delta_{C}(A) = A^{n} + d_{1} A^{n-1} + \ldots + d_{n-1} A + d_{n}I.
$$

\smallskip

\noindent
By eliminating $y_{2}, \ldots, y_{n}$ from the previous system,
we obtain the operator equation$:$
\begin{equation}
\begin{array}{rcl}
\label{y-Op-jednacina}
\Delta_{C}(A)(y_{1}) &\!\!=\!\!& {\Big (}
A^{n-1}(\psi_{1})
{\Big )}                                                                                                \\[0.75 ex]
&\!\!+\!\!& {\Big (} A^{n-2}(\psi_{2}) + d_{1} A^{n-2}(\psi_{1})
{\Big )}                                                                                                \\[0.75 ex]
&\!\!+\!\!& {\Big (} A^{n-3}(\psi_{3}) + d_{1} A^{n-3}(\psi_{2}) +
d_{2} A^{n-3}(\psi_{1})
{\Big )}                                                                                                \\[-0.5 ex]
&\!\!\vdots\!\!&                                                                                        \\[-0.5 ex]
&\!\!+\!\!& {\Big (} \psi_{n} + d_{1} \psi_{n-1} + d_{2}
\psi_{n-2} + \ldots + d_{n-2} \psi_{2} + d_{n-1} \psi_{1}
{\Big )}.
\end{array}
\end{equation}
Using formulae (\ref{C-Sistem-1}), (\ref{y-Op-jednacina}) and Lemma \ref{Lema3} we complete the proof. $ \; \Box$

\vspace{5pt}
\begin{remark}{\rm
It is important to say that the reduction from the previous theorem is reversible,
i.e. reduced system {\rm (\ref{PRS 1})} implies system {\rm (\ref{B-Sistem-1})}.}
\end{remark}

\vspace{5pt}
\begin{remark}{\rm
System {\rm (\ref{PRS 1})} is called {\bf partially reduced}
because it has only one variable separated from the others by a
linear $n^{th}$ order $A$-operator equation.}
\end{remark}

\medskip

We now continue with the general case of partial reduction.

\smallskip

\begin{theorem}
\label{Teorema_parc_red_2}
Assume that the rational canonical form of the matrix $B$ is the block diagonal matrix
$$
C = \mbox{\rm diag}{\big (}C_{1}, C_{2}, \ldots, C_{k}{\big )} \;\;\; \mbox{for some}\; k, \; 2 \leq k \leq n,
$$
where $C_{i}$ are the companion matrices of the monic polynomials
$\Delta_{C_{i}}(\lambda)\!=\!
\lambda^{n_{i}} + d_{i,1}\lambda^{n_{i}-1} + \ldots + d_{i,n_{i}-1} \lambda + d_{i,n_{i}}$
of degree $n_{i} \geq 1$ and $\Delta_{C_{i}} \;|\; \Delta_{C_{i+1}}$ for all $i$, $1 \leq i < k$.

\smallskip

\noindent
Let $\ell_{1} \!=\! 0$ and let $\ell_{i} \!=\! \sum_{j=1}^{i-1}{n_{j}}$ for $i$, $2 \leq\! i\! \leq k$.
Then linear system of first order $A$-operator equations {\rm (\ref{B-Sistem-1})}
can be transformed into the conjunction of the partially reduced systems$:$
\begin{equation}
\label{wedge 1}
\hspace{180pt}
\bigwedge\limits_{i=1}^{k}{
(\mbox{$\cal R$}_{C_{i}})},
\end{equation}

\break

\noindent
where every subsystem $(\mbox{$\cal R$}_{C_{i}})$ is of the form$:$
$$
\left\{
\begin{array}{ccl}
\Delta_{C_{i}}(A)(y_{\ell_{i}+1}) &\!=\!& \displaystyle\sum\limits_{k=1}^{n_{i}}{\!(-1)^{k+1}
\mbox{\large $\delta$}^{1}_{k}(C_{i}; A^{n_i-k}(\psi_{\ell_{i}+1}), \ldots, A^{n_i-k}(\psi_{\ell_{i}+n_{i}}))\!} \\[0.5 ex]
y_{\ell_{i}+2}            &\!=\!&  A(y_{\ell_{i}+1}) - \psi_{\ell_{i}+1}                                         \\[0.0 ex]
y_{\ell_{i}+3}            &\!=\!&  A(y_{\ell_{i}+2}) - \psi_{\ell_{i}+2}                                         \\[-0.5 ex]
                          & \vdots &                                                                             \\[-0.5 ex]
y_{\ell_{i}+n_{i}}        &\!=\!&  A(y_{\ell_{i}+n_{i}-1}) - \psi_{\ell_{i}+n_{i}-1}
\end{array}
\right\}.
$$
Columns $\vec{y}\!=\![y_{1}\;\ldots\;y_{n}]^{T}$ and $\vec{\psi} \!=\! [\psi_{1}\;\ldots\;\psi_{n}]^{T}$
are determined by $\vec{y} \!=\! P^{-1} \vec{x}$ and $\vec{\psi} \!=\! P^{-1} \vec{\varphi}$
for a regular matrix $P$ such that $C \!=\! P^{-1} B P$.
\end{theorem}

\vspace{5pt}
\noindent
{\bf Proof.}
According to Lemma \ref{Lema1}, system (\ref{B-Sistem-2}) is equivalent to the following system$:$
\begin{equation}
\label{wedge_2}
\hspace{180pt}
\bigwedge\limits_{i=1}^{k}{{\bigg (}\vec{A}(\vec{y}_{i}) = C_{i} \vec{y}_{i} + \vec{\psi}_{i}{\bigg )}},
\end{equation}
where
$\vec{y}_{i} \!=\![y_{\ell_{i} + 1}\;\ldots\;y_{\ell_{i} + n_{i}}]^{T}$
is a column of unknowns and
$\vec{\psi}_{i} \!=\![\psi_{\ell_{i} + 1}\;\ldots\;\psi_{\ell_{i} + n_{i}}]^{T}$
is the free column of the linear system
$\vec{A}(\vec{y}_{i}) = C_{i} \vec{y}_{i} + \vec{\psi}_{i}$. We now consider the linear operator$:$
$$
\hspace{-200pt}
\Delta_{C_{i}}(A)=
A^{n_{i}} + d_{i,1}A^{n_{i}-1} + \ldots + d_{i,n_{i}-1} A + d_{i,n_{i}}I .
$$
Therefore, by Theorem \ref{Teorema_parc_red_1} every subsystem
$\vec{A}(\vec{y}_{i}) = C_{i} \vec{y}_{i} + \vec{\psi}_{i}$
is equivalent to the partially reduced system $(\mbox{$\cal R$}_{C_{i}})$.
We now obtain the equivalence between system (\ref{wedge_2}) and the system
$\bigwedge\limits_{i=1}^{k}{(\mbox{$\cal R$}_{C_{i}})}. \; \Box $

\begin{remark}{\rm
System $\bigwedge\limits_{i=1}^{k}{(\mbox{$\cal R$}_{C_{i}})}$ is equivalent to the initial
system {\rm (\ref{B-Sistem-1})}. Every subsystem $(\mbox{$\cal R$}_{C_{i}})$ consists of
a higher order linear $A$-operator equation, in only one unknown,
and first order linear $A$-operator equations, in exactly two variables.}
\end{remark}

\vspace{5pt}
\begin{remark}{\rm
The initial system has a solution if and only if every higher order linear
$A$-operator equation of the partially reduced system has a solution.}
\end{remark}

\vspace{5pt}
\begin{remark}{\rm
The above mentioned reduction formulae can be applied to two special types of non-homogeneous systems, to systems of
differential and systems of difference equations with constant coefficients.
The usual way of proving the reduction formulae for a non-homogeneous system of differential equations with constant
coefficients is to use the Jordan canonical form. However, this method is somewhat complicated, see {\bf [3, {\rm pp. 166-178}]}.
We give a better and more efficient alternative formulae of reduction for an arbitrary \mbox{$A$-operator}.
Our method involves the rational canonical instead of Jordan form which makes reduction process much more easier
to handle and enables us to find general solutions.}
\end{remark}

\section{Examples}

The reduction formulae in special cases of $n=2$ and $n=3$ will be given as examples.

\vspace{5pt}
\begin{example}
{\rm
Consider the linear system of the first order $A$-operator equations, in unknowns $x_{1}, x_{2}$$:$
\begin{equation}
\label{Primer2}
\hspace{90pt}
\left\{
\begin{array}{c}
A(x_{1}) = b_{11} x_{1} + b_{12} x_{2} + \varphi_{1} \\[0.5 ex]
A(x_{2}) = b_{21} x_{1} + b_{22} x_{2} + \varphi_{2}
\end{array}
\right\}.
\end{equation}
Depending on the block structure of the rational canonical form $C$ of the matrix $B$,
partially reduced system of system (\ref{Primer2}) has two possible forms.

\pagebreak

{\boldmath $1^{0}$}
The rational canonical form $C$ of the matrix $B$ has one block.
Then the partially reduced system has

\hspace{10pt}
the form$:$

\vspace{-5pt}
\hspace{100pt}
$
\left\{
\begin{array}{l}
\Delta_{B}(A)(y_1) =
\mbox{\large $\delta$}^{1}_{1}(C; A(\psi_{1}),A(\psi_{2}))
\!-\!
\mbox{\large $\delta$}^{1}_{2}(C; \psi_{1}, \psi_{2}) \\ [0.5 ex]
y_2 = A(y_1) - \psi_{1}
\end{array}
\right\}.
$

\vspace{5pt}
{\boldmath $2^{0}$}
The rational canonical form $C$ of the matrix $B$ is diagonal matrix $C=\mbox{diag}(C_{1}, C_{1})$.
Then system (\ref{Primer2})

\hspace{10pt}
is equivalent to the partially reduced system$:$

\hspace{100pt}
$
\left\{
\begin{array}{l}
\Delta_{C_1}(A)(y_1) = \psi_{1} \\[0.5 ex]
\Delta_{C_1}(A)(y_2) = \psi_{2}
\end{array}
\right\}.
$
}
\end{example}

\vspace{5pt}
\begin{example}{\rm
We now deal with the system of first order $A$-operator equations in three
unknowns $x_{1}, x_{2}, x_{3}$$:$
\begin{equation}
\label{Primer3}
\hspace{90pt}
\left\{
\begin{array}{c}
A(x_1) = b_{11} x_{1} + b_{12} x_{2} + b_{13} x_{3} + \varphi_{1}    \\[0.5 ex]
A(x_2) = b_{21} x_{1} + b_{22} x_{2} + b_{23} x_{3} + \varphi_{2}    \\[0.5 ex]
A(x_3) = b_{31} x_{1} + b_{32} x_{2} + b_{33} x_{3} + \varphi_{3}
\end{array}
\right\}.
\end{equation}

\noindent
Partially reduced system (\ref{wedge 1}) of system (\ref{Primer3}) depends on the block structure
of the rational canonical form~$C$ of the matrix $B$. The following cases are possible.

\vspace{5pt}
{\boldmath $1^{0}$}
The rational canonical form $C$ of the matrix $B$ has one block; then (\ref{wedge 1}) has the form

\vspace{5pt}
\hspace{100pt}
$
\left\{
\begin{array}{l}
\Delta_{B}(A)(y_1) =
\displaystyle\sum\limits_{k=1}^{3}{\!(-1)^{k+1}
\mbox{\large $\delta$}^{1}_{k}(C; A^{3-k}(\psi_{1}), A^{3-k}(\psi_{2}), A^{3-k}(\psi_{3}))\!} \\ [0.5 ex]
y_2 = A(y_1) - \psi_{1}                                                                       \\ [0.5 ex]
y_3 = A(y_2) - \psi_{2}
\end{array}
\right\}.
$

\medskip

{\boldmath $2^{0}$}
The rational canonical form $C$ of the matrix $B$ has two blocks $C=\mbox{diag}(C_{1}, C_{2})$,
where~$C_{1}$~and~$C_{2}$

\hspace{11pt}
are matrices of order 1 and 2 respectively; then (\ref{wedge 1}) has the form

\vspace{5pt}
\hspace{100pt}
$
\left\{
\begin{array}{l}
\Delta_{C_{1}}(A)(y_1) = \psi_{1}                                            \\ [0.5 ex]
\Delta_{C_{2}}(A)(y_2) = \mbox{\large $\delta$}^{1}_{1}(C_{2}; A(\psi_{2}), A(\psi_{3}))
- \mbox{\large $\delta$}^{1}_{2}(C_{2}; \psi_{2}, \psi_{3})                  \\ [0.5 ex]
y_3 = A(y_2) - \psi_{2}
\end{array}
\right\}.
$

\medskip

{\boldmath $3^{0}$}
The rational canonical form $C$ of the matrix $B$ has three blocks equal to $C_{1}$;
then (\ref{wedge 1}) has the form

\vspace{5pt}
\hspace{100pt}
$
\left\{
\begin{array}{c}
\Delta_{C_{1}}(A)(y_1) = \psi_{1} \\
\Delta_{C_{1}}(A)(y_2) = \psi_{2} \\
\Delta_{C_{1}}(A)(y_3) = \psi_{3}
\end{array}
\right\}.
$
}
\end{example}

\bigskip

\smallskip
\noindent
{\bf Acknowledgment.} The authors would like to thank the reviewers for their valuable suggestions and comments.
Research is partially supported by the MNTRS, Serbia, Grant No. 144020.

\vspace{2cc}
\begin{center}{\bf REFERENCES}
\end{center}

\vspace{0.8cc}
\newcounter{ref}
\begin{list}{\bf [\arabic{ref}]}{\usecounter{ref} \leftmargin 4mm
\itemsep -1mm}

\item{F.M. Goodman, {\em Algebra: abstract and concrete}, Semisimple Press, Iowa City, 1998.}\\[-1.00 ex]

\item{{\cirilrm F.R. Gantmaher}, {\cirilit Teori{\ja } matric}, {\cirilrm izdanie \ch etvertoe}, {\cirilrm Nauka}, {\cirilrm Moskva}, 1988.}\\[-1.00 ex]

\item{{\cirilrm I.G. Petrovski\ijj}, {\cirilit Lekcii po teorii ob\ii knovenn\ii h differencial\mek n\ii h uravneni\ijj},
{\cirilrm izdanie sed\mek moe}, {\cirilrm Izdatel\mek stvo moskovskogo universiteta}, {\cirilrm Moskva}, 1984.}\\[-1.00 ex]

\item{J.C. Butcher, {\em A Mathematical Miniature}, Mathematical Miniatures and Apologies, NZMS Newsletter, No.~68, December 1996.
({\sf http:/\mbox{$\!$}/www.math.auckland.ac.nz/ \mbox{$\!\!\sim$}butcher/miniature/miniatureAndApology.html})}\\[-1.00 ex]

\item{J.C. Butcher, P. Chartier, {\em The effective order of singly implicit Runge Kuta methods}, Numerical Algorithms 20 (4) (1999)
269-284.}\\[-1.00 ex]

\item{J.C. Butcher, W.M. Wright, {\em Applications of doubly companion matrices}, Applied Numerical Mathematics 56 (3) (2006) 358-373.}\\[-1.00 ex]

\item{J.C. Butcher, {\em General linear methods for ordinary differential equations}, Mathematics and Computers in Simulation 79 (6) (2009) 1834-1845.}

\end{list}

\end{document}